\documentclass[reqno, 12pt]{amsart}
\usepackage{epsfig,amsmath,amsfonts,latexsym}
\usepackage{amsthm}
\usepackage[abs]{overpic}
\usepackage[usenames,dvipsnames]{color}
\usepackage{palatino}
\usepackage[hyperfootnotes=false]{hyperref}
\usepackage{caption}
\usepackage{dsfont}
\usepackage{mathtools}
\usepackage{cite}
\usepackage{tikz-cd}
\usepackage{amssymb}
\usepackage{comment}
\usepackage{epigraph}
\hypersetup{
  colorlinks,
  citecolor=Blue,
  linkcolor=Black,
  urlcolor=Green}
  
\theoremstyle{plain}

\newtheorem{theorem}{Theorem}

\newtheorem{dummy}{anything}[section]

\newtheorem{proposition}[dummy]{Proposition}

\theoremstyle{definition}
\newtheorem{definition}[dummy]{Definition}

\newtheorem{example}[dummy]{Example}

\theoremstyle{remark}
\newcommand{\R}{\mathbb{R}}

\def\R{\mathbb{R}}

\title{Pseudo-holomorphic dynamics in the three-body problem: fillability and convexity} 

\author{Agustin Moreno}

\address[A.\ Moreno]{School of Mathematics \\ Institute for Advanced Study \\ 1 Einstein Dr, 08540, Princeton NJ  \\ USA \& Mathematisches Institut \\ Universit\"at Heidelberg\\ Im Neuenheimer Feld 205\\ 69120 Heidelberg \\ Germany }

\email{agustin.moreno2191@gmail.com}

\date{}

\begin{document}

\maketitle

\begin{abstract} In this article, we extend the methods from \cite{M}, where the five dimensional analogue of the three dimensional finite energy foliations introduced by Hofer--Wysocki--Zehnder was identified, to the case where there the underlying (IP) contact $5$-fold admits a $6$-dimensional (IP) symplectic filling. We show that the filling induces a moduli space of pseudo-holomorphic curves which is \emph{itself} a symplectic filling of the standard $3$-sphere, and hence symplectomorphic to the $4$-dimensional ball. We further show that whenever the contact form on the $5$-fold is strictly convex, then this moduli space is a strictly convex domain, so that the induced Reeb dynamics at the boundary is in particular dynamically convex. For the circular restricted three-body problem, this implies that whenever the spatial dynamics (near one of the heavy masses) is strictly convex, the holomorphic shadow of \cite{M} is dynamically convex; this is shown to hold for near-integrable cases close to the Kepler problem, and with mass ratio either zero or sufficiently close to 1.
\end{abstract}

\section{Introduction}

This article is a follow-up of \cite{M}, where the theory of pseudo-holomorphic curves was leveraged in order to obtain $5$-dimensional finite energy foliations for the spatial (circular) restricted three-body problem, which are higher-dimensional analogues of the foliations introduced by Hofer--Wysocki--Zehnder for the study of three-dimensional Reeb flows. They were used to associate, to the spatial dynamics of the three-body problem (near one of the heavy masses), a Reeb flow on the standard three-sphere which is in some sense a lower-dimensional ``shadow'' of the original dynamics. In this article, we build on these methods, to deal with the case where the $5$-dimensional contact manifold admits a $6$-dimensional filling which is compatible with the underlying \emph{iterated planar} (IP) structure. This structure is indeed present in the spatial restricted three-body problem, as follows from \cite{MvK} (see \cite{M}). We will then be interested in properties of the associated shadow; in particular, we will give natural conditions for its dynamical/strict convexity.  

\medskip

\medskip

\textbf{Setup.} As in \cite{M}, we consider an iterated planar (IP) contact $5$-fold $(M,\xi)$. In particular, we have a concrete open book decomposition $\theta_M: M \backslash B \rightarrow S^1$ on $M$, where $B$ is the binding of the open book, which supports $\xi$ in the sense of Giroux. That is, there is a contact form $\alpha$ for $\xi$, a \emph{Giroux form}, such that $\alpha_B=\alpha\vert_B$ is contact, and $d\alpha$ is positively symplectic on the fibers of $\theta_M$. This is  equivalent to the Reeb flow of $\alpha$ having $B$ as an invariant subset, and being positively transverse to each fiber. We use the abstract notation $(M,\xi)=\mathbf{OB}(P,\phi)$, where $P$ denotes the abstract page (the closure of the typical fiber of $\theta_M$) with $\partial P=B$, and $\phi$ is the symplectic monodromy. The IP structure means further that $P$ admits the structure of a $4$-dimensional Lefschetz fibration over $\mathbb{D}^2$ whose fibers are surfaces of genus zero. We write $P=\mathbf{LF}(F,\phi_F)$, where $\phi_F$ is the monodromy of the Lefschetz fibration on $P$, and $F$ is its typical fiber, having genus zero and as potentially several boundary components. The binding is a planar contact $3$-fold, i.e.\ we have $(B,\xi_B)=\mathbf{OB}(F,\phi_F)$.

As in \cite{M}, we assume that the concrete open book on $M$ is adapted to the Reeb dynamics of a \emph{fixed} contact form $\alpha$, whose dynamics we wish to study, which we assume is a \emph{IP Giroux form}. This means that we have a further concrete planar open book $\theta_B: B\backslash L\rightarrow S^1$ on the $3$-manifold $B=\mathbf{OB}(F,\phi_F)$ of the given abstract type, which is adapted to the Reeb dynamics of $\alpha_B$. Then $L=\partial F$ is a link in $B$, the binding of the open book for $B$ (i.e.\ the "binding of the binding''), consisting of Reeb orbits for $\alpha_B$. 

For $F$ a genus zero surface, let $\mathbf{Reeb}(F,\phi_F)$ be the collection of contact forms whose flow is adapted to some concrete planar open book $\theta_B: B\backslash L \rightarrow S^1$ on $B$, of abstract form $B=\mathbf{OB}(F,\phi_F)$. Iteratively, we let $\mathbf{Reeb}(\mathbf{LF}(F,\phi_F),\phi)$ be the collection of contact forms with flow adapted to some concrete IP open book $\theta_M: M\backslash B \rightarrow S^1$ on manifold $M$, of abstract form $M=\mathbf{OB}(\mathbf{LF}(F,\phi_F),\phi)$, whose restriction to the binding $B=\mathbf{OB}(F,\phi_F)$ belongs to $\mathbf{Reeb}(F,\phi_F)$.

In this setup, \cite[Theorem A]{M} gives a moduli space $\overline{\mathcal{M}}=\mathbb R\times \overline{\mathcal{M}}^q$ of planar pseudo-holomorphic curves in $\mathbb R\times M$. The projections of these curves to $M$ are the leaves of a foliation $\overline{\mathcal{M}}^q$ of $M\backslash L$, forming the fibers of a concrete Lefschetz fibration on each page $P_\varphi=\overline{\theta_M^{-1}(\varphi)}$ of the open book on $M$ (of abstract type $\mathbf{LF}(F,\phi_F)$), which are symplectic with respect to $\omega_\varphi=d\alpha\vert_{P_\varphi}$. We then write $(P_\varphi,\omega_\varphi)=\mathbf{LF}(F,\phi_F)$ to indicate that the symplectic structure is compatible with the Lefschetz fibration, i.e.\ the symplectic form is positive on the fibers. Along the binding $B$, the leaves of $\overline{\mathcal{M}}^q$ form the pages $\mathcal{M}^q_B\cong S^1$ of the open book on $B$. Moreover, the leaf space is $\overline{\mathcal{M}}^q\cong S^3=\mathbf{OB}(\mathbb D^2,\mathds 1)$, endowed with a concrete and trivial open book $\theta_\mathcal{M}: \overline{\mathcal{M}}^q\backslash \mathcal{M}^q_B \rightarrow S^1$, with binding $\mathcal{M}^q_B\cong S^1$, and whose page $\mathbb D^2_\varphi=\overline{\theta_\mathcal{M}^{-1}(\varphi)}$ correspond to the base of the Lefschetz fibration of $P_\varphi$.

In this paper, we will extend the above setup of \cite{M}, and further assume that we have a $6$-dimensional exact symplectic filling $(W,\omega=d\lambda)=\mathbf{LF}(P,\phi)$ of $M$, whose symplectic form is adapted to a concrete Lefschetz fibration $\pi_W: W \rightarrow \mathbb D^2$ of the given abstract type, inducing the IP open book at the boundary $M=\partial W$. We also assume that there is a symplectic submanifold $(P_B,\omega_B=\omega\vert_{P_B})\subset W$ such that $W=\mathbf{LF}(F,\phi_F)$ abstractly, which gives a filling of $(B,\xi_B)=\partial P_B$, and which is a regular fiber of $\pi_W$ (so in particular $P_B\cong P$). We call the data $((W,\omega)=\mathbf{LF}(P,\phi), (P_B,\omega_B)=\mathbf{LF}(F,\phi_F))$ an \emph{IP filling} of the IP contact manifold $((M,\xi)=\mathbf{OB}(P,\phi),(B,\xi_B)=\mathbf{OB}(F,\phi_F))$. Note that we do not fix concrete Lefschetz fibrations for $P$ nor for $P_B$, as we will construct them.

\smallskip

\begin{example}
    The case relevant for the restricted three-body problem, for the dynamics near one of the primaries, is $W=D^*S^3=\mathbf{LF}(D^*S^2,\tau^2)$ (where $\tau$ is the Dehn-Seidel twist), which is a filling of the Moser-regularized spatial problem $S^*S^3=\mathbf{OB}(D^*S^2,\tau^2)$; and $P_B=D^*S^2=\mathbf{LF}(D^*S^1,\tau_0^2)$ (where $\tau_0$ is the Dehn twist), a filling of the Moser-regularized planar problem $S^*S^2=\mathbf{OB}(D^*S^1,\tau_0^2)$. The concrete open book for the spatial problem was found by the author and Otto van Koert in \cite{MvK}, while that for the planar problem, was found in \cite{HSW} whenever the planar dynamics is dynamically convex.
\end{example}

\textbf{Statement of results.} The first result of this paper extends the foliation $\overline{\mathcal{M}}^q\cong S^3$ to a foliation $\mathcal{\overline{M}}^W\cong \mathbb D^4$ on $W$, such that $\partial\mathcal{\overline{M}}^W=\overline{\mathcal{M}}^q$, and inducing concrete Lefschetz fibrations for $P$ and $P_B$. 

\begin{theorem}[IP foliation]\label{thm:IPfoliation} There is a foliation $\overline{\mathcal{M}}^W$ of $W\backslash L$, consisting of immersed symplectic surfaces whose boundary is $L$. We have $\overline{\mathcal{M}}^W\cong \mathbb D^4=\mathbf{LF}(\mathbb D^2,\mathds 1)$, and its boundary is naturally identified with $\partial \overline{\mathcal{M}}^W=\overline{\mathcal{M}}^q\cong S^3=\mathbf{OB}(\mathbb D^2,\mathds 1)$. The curves in $\overline{\mathcal{M}}^W$ foliate each regular fiber of the concrete Lefschetz fibration on $W=\mathbf{LF}(P,\phi)$, along which they themselves form the fibers of a concrete Lefschetz fibration of abstract type $\mathbf{LF}(F,\phi_F)$. In particular, there is a subfoliation $\overline{\mathcal{M}}_B\subset \overline{\mathcal{M}}^W$ foliating $P_B$, with $\overline{\mathcal{M}}_B\cong \mathbb D^2$, which is the typical fiber of a concrete trivial Lefschetz fibration $\pi_\mathcal{M}:\overline{\mathcal{M}}^W \rightarrow \mathbb D^2$, and whose elements consist of fibers of a concrete Lefschetz fibration $\theta_B:P_B\rightarrow \mathbb D^2$ of abstract type $\mathbf{LF}(F,\phi_F)$.
    
\end{theorem}

We call a foliation as above, an \emph{IP foliation} of the IP filling. We denote $\widehat{\mathcal M}=\overline{\mathcal M}^W\cup\overline{\mathcal M}$, which is a foliation of the Liouville completion $\widehat W=W\cup [0,+\infty)\times M$ by curves which are asymptotically cylindrical, with asymptotics corresponding to the link $L\subset B \subset M$.

\medskip

\textbf{Contact and symplectic structures on moduli.} From \cite[Theorem B]{M}, we know that $\overline{\mathcal{M}}^q$ carries a natural contact structure $\xi_\mathcal{M}$ which is supported by the trivial open book on $S^3$, and therefore isomorphic to the standard contact structure. Moreover, the Giroux form $\alpha$ on $M$ induces a tautological contact form $\alpha_\mathcal{M}$ on $\overline{\mathcal{M}}^q$ by leaf-wise integration (called the \emph{holomorphic shadow}), whose Reeb flow is adapted to the trivial open book on $\overline{\mathcal{M}}^q$. The next result in this paper extends the above geometric structures to $\overline{\mathcal{M}}^W$. 

\begin{theorem}[contact and symplectic structures on moduli]\label{thm:symplectic_structure} The moduli space $\overline{\mathcal{M}}^W\cong \mathbb D^4$ carries a tautological symplectic form $\omega_\mathcal{M}$ which gives a symplectic filling of its strict contact boundary $(\overline{\mathcal{M}}^q,\alpha_\mathcal{M})$, and is therefore symplectomorphic to a domain in the standard symplectic $\mathbb C^2$. The symplectic structure extends naturally to $\widehat{\mathcal{M}}$, which is the liouville completion of $\overline{\mathcal{M}}^W$, and symplectomorphic to the standard $\mathbb C^2$.
    
\end{theorem}

Here, recall that a celebrated result of Gromov implies that a (minimal) filling of the standard $S^3$ is symplectomorphic to the standard $4$-ball.

\medskip

\textbf{Convexity.} In the following, we will investigate when the holomorphic shadow defined in \cite{M} (i.e.\ the contact form at the boundary of $\overline{\mathcal{M}}^W$) is given by a dynamically convex contact form on $S^3$. We will arrive at this notion via the stronger notion of \emph{strict convexity}, which we now recall.

\begin{definition}[strict convexity]\label{def:convexity}
    Let $(M,\alpha)$ be a strict contact-type hypersurface in a symplectic manifold $(W,\omega)$. We say that $(M,\alpha)$ (or its Reeb flow) is \emph{strictly convex} if $M=H^{-1}(0)$, for a $C^2$-Hamiltonian $H:(-\epsilon,\epsilon)\times M\subset W\rightarrow \mathbb R$ defined near $M$, satisfying $\nabla^2H\vert_M$ is positive definite. Here $\nabla^2H\vert_M=\nabla (dH)\vert_M: TM\rightarrow T^*M$ is the Hessian of $H$, given by $\nabla^2H\vert_M= \nabla(dH)\vert_{T^*M}$.
\end{definition}

Note that $dH\vert_M=0$ in the above, so that $\nabla^2H=\nabla (dH)$ is well-defined along $M$. The connection used is the Levi-Civita connection of some ambient metric on $W$\footnote{Note that if $\nabla^M$ is the Levi-Civita connection of the metric induced on $M$, we have $$\nabla_XY= \nabla^M_XY+\mathrm{II}(X,Y),$$ for $X,Y\in TM$, where $\mathrm{II}$ is the second fundamental form of $M$. Therefore
$$
\nabla^M_XY=\pi \nabla_XY,
$$
where $\pi$ is the orthogonal projection to $TM$.}. If $M$ is three-dimensional, we also have the following notion of dynamical convexity.

\begin{definition}[Dynamical convexity]
    Let $(M,\alpha)$ be a strict contact $3$-manifold with $c_1(\xi)\vert_{\pi_2(M)}=0$. Then $(M,\alpha)$ is dynamically convex if the Conley-Zehnder index of every contractible Reeb orbit $\gamma$ (with respect to a disk capping trivialization) satisfies $\mu_{CZ}(\gamma)\geq 3$.
\end{definition}

The two notions are related in one direction, i.e.\ strict convexity implies dynamical convexity \cite{HWZ}. Note that the converse is not true \cite{CE}. In our setup, we have the following.

\begin{theorem}[Convexity] \label{thm:convexity} Assume that the strict IP $5$-fold $(M,\alpha)$ is a strictly convex contact-type hypersurface in the IP filling $(W,\omega)$. Then the holomorphic shadow $(\overline{\mathcal{M}}^q,\alpha_\mathcal{M})$ is a strictly convex contact-type hypersurface in $(\overline{\mathcal M}^W,\omega_\mathcal{M})$. In particular, it is a dynamically convex contact $3$-sphere.
    
\end{theorem}

\textbf{Restricted three-body problem.} Recall that the (circular) restricted 3BP concerns the motion of a massless particle in $\mathbb R^3$ under the gravitational influence of two heavy primaries with mass $\mu$ and $1-\mu$, where $\mu\in(0,1)$, which move in circles around their center of mass. In a rotating frame centered at the latter, we have $m=(\mu-1,0,0)$ (the Moon), and $e=(\mu,0,0)$ (the Earth). The problem is given by the Hamiltonian 
$H:(T^*\R^3\setminus \{  m, e\},d p \wedge d  q)\rightarrow \mathbb R$
$$
H( q, p)=\frac{1}{2}\Vert  p\Vert^2 - \frac{\mu}{\Vert  q- m\Vert } - \frac{1-\mu}{\Vert  q- e\Vert } +p_1q_2-p_2q_1. 
$$
The planar case of this problem is obtained by setting $q_3=p_3=0$; the Jacobi constant $c$ is the value of $H$. The Hamiltonian flow of this dynamical system has singularities caused by two-body collisions. However, the resulting flow can be extended across these singularities using e.g.\ Moser regularization \cite{moser}. If $L_1$ is the Lagrange critical points with smallest energy value, then if $c<H(L_1)$, the Moser-regularized energy level set near the Earth or Moon is topologically $S^*S^3$, with the invariant planar problem given by $S^*S^2$. These components of the level set are contact-type \cite{AFvKP}, carrying the corresponding standard contact structures. Moreover, they fit into an IP structure $$((S^*S^3,\xi_{std})=\mathbf{OB}(D^*S^2,\tau^2), (S^*S^2,\xi_{std})=\mathbf{OB}(D^*S^1,\tau_0^2)),$$ with a natural IP filling $$((D^*S^3,\omega_{std})=\mathbf{LF}(D^*S^2,\tau^2), (D^*S^2,\omega_{std})=\mathbf{OB}(D^*S^1,\tau_0^2)).$$ Moreover, we have a \emph{concrete} open book adapted to the dynamics of the spatial problem, as follows from \cite{MvK}. And if the planar dynamics is dynamically convex, we also have a concrete open book adapted to the planar problem \cite{HSW}. This was the starting point for the results in \cite{M} which we are now extending. 

In the limit $c\rightarrow -\infty$, Moser regularization recovers the Kepler problem, which is given by the geodesic flow of the round metric on $S^3$, having the planar problem $S^2$ as a totally geodesic submanifold. For this limit case, we have the following elementary result, which follows by a computation which was communicated to the author by Otto van Koert.

\begin{proposition}[Kepler problem]\label{prop:Kepler_problem_convexity}
    After Moser regularization, the negative energy Kepler problem is a strictly convex contact-type hypersurface in $(T^*S^3,\omega_{std})$. 
\end{proposition}

It was observed in \cite{M} that the shadow of the Kepler problem is the Hopf flow on $S^3$, which is indeed a strictly convex Reeb flow. As convexity is a $C^2$ open condition, convexity of the spatial problem near one of the heavy masses will still be true for $c\ll 0$. Moreover, for $c\ll 0$ and mass ratio $\mu$ either zero or sufficiently close to $1$, the planar problem is still dynamically convex (it is actually strictly convex in the latter case) \cite{AFFHvK,AFFvK}, so that we have the concrete planar open book provided by \cite{HSW}, and then can apply the results in this paper --in particular, there is a well-defined shadow dynamics for this range of parameters--. The following perturbative result is then a corollary of Theorem \ref{thm:convexity}.

\begin{theorem}[Near integrable convexity] If $c\ll 0$, then the spatial problem near one of the heavy primaries is given by a strictly convex Reeb flow in $S^*S^3$. If moreover $\mu$ is either zero or sufficiently close to $1$, then its holomorphic shadow is a strictly convex Reeb flow in $S^3$, and in particular, dynamically convex.
    
\end{theorem}

\textbf{Acknowledgements.} The author is grateful to Umberto Hryniewicz for productive conversations and interest in the project. The author received support by the National Science Foundation under Grant No.\ DMS-1926686, and is currently supported by the Sonderforschungsbereich TRR 191 Symplectic Structures in Geometry, Algebra and Dynamics, funded by the DFG (Projektnummer 281071066 – TRR 191), and also by the DFG under Germany's Excellence Strategy EXC 2181/1 - 390900948 (the Heidelberg STRUCTURES Excellence Cluster). 

\section{Convexity}

We first present the proof of Proposition \ref{prop:Kepler_problem_convexity}.

\medskip

\textbf{Round geodesic flow.} We now show that the regularized spatial Kepler problem in negative energy is strictly convex. Indeed, this is simply the geodesic flow of the round metric on $S^3$, so it suffices to show this for the round geodesic flow on $S^n$ for all $n\geq 1$. The following computations were communicated to us by Otto van Koert.

We write
$$
T^*S^n=\left\{(\xi;\eta)\in \mathbb R^{n+1}\oplus \mathbb R^{n+1}: \Vert\xi \Vert=1, \langle \xi,\eta\rangle = 0\right\},
$$
endowed with the symplectic form $\omega=d\xi \wedge d\eta\vert_{T^*S^n}$, the round metric, and the standard almost complex structure $J$ with $J \partial_{\xi}=\partial_\eta$. The Hamiltonian is then
$$
Q: T^*S^n\rightarrow \mathbb R,
$$
$$
Q(\xi;\eta)= \frac{1}{2}\Vert \eta \Vert^2\vert_{T^*S^n}. 
$$
We now compute the Hessian of $Q$ along $M:=S^*S^n=H^{-1}(1/2)$, using that 
$$
\nabla X_Q=\nabla (J\nabla H)=J\nabla^2H,
$$
so that
$$
\nabla^2H=-J\nabla X_Q.
$$
Here we have used that $\nabla J=0$.

To simplify the computations, we first note that $SO(n+1)$ acts transitively on $S^*S^n$, so we may assume that we evaluate the Hessian at the point
$$
p=(1,0,\ldots,0;0,1,0,\ldots ,0) =(\xi;\eta)\in M.
$$

Set $W:= T^*S^n \subset T^* \mathbb R^{n+1}$. Then 
$$
TW=\ker(\xi \cdot d\xi)\cap \ker (\xi \cdot d\eta + \eta \cdot d\xi).
$$
We take the basis of $T_pW=\ker(d\xi_0)\cap \ker(d\eta_0+d\xi_1)$ given by
$$
U_1:=\partial_{\xi_1} -\partial_{\eta_0},\; V_1:=   \partial_{\eta_1},
$$
$$
U_j:=\partial_{\xi_j},\; V_j:=\partial_{\eta_j},\mbox{ for } j\geq 2.
$$

We observe that these form a symplectic basis. We then have
$$
X_Q =\eta \cdot \partial_\xi -\Vert \eta\Vert^2 \xi\cdot \partial_\eta
$$

We denote the Levi-Civita connection on $\mathbb R^{2n+2}$ by $\nabla$, and that of the metric induced on $W$, by $\nabla^W$. Let $\pi: T(\mathbb R^{2n+2})\vert_M\rightarrow TM$ denote the orthogonal projection.

For $j\geq 2$, along $M$, and evaluated at $p$ (which has $\eta_j=0$), we have
$$
\nabla^W_{U_j} X_Q\vert_p=\pi\nabla_{U_j} X_Q\vert_p=-\partial_{\eta_j}=-V_j,
$$
$$
\nabla^W_{V_j} X_Q\vert_p=\pi\nabla_{V_j} X_Q\vert_p=\partial_{\xi_j}-2\eta_j (\xi \cdot \partial_\eta)\vert_p=U_j.
$$
For $j=1$, we have
$$
\nabla^W_{U_1} X_Q\vert_p=\pi\nabla_{U_1} X_Q\vert_p=\pi(-\partial_{\eta_1}-\partial_{\xi_0})=-V_1.
$$
$$
\nabla^W_{V_1} X_Q\vert_p=\pi\nabla_{V_1} X_Q\vert_p=\pi(\partial_{\xi_1}-2\partial_{\eta_0}
)=\left\langle \partial_{\xi_1}-2\partial_{\eta_0}, \frac{1}{\sqrt{2}}U_1\right\rangle\frac{1}{\sqrt{2}}U_1=\frac{3}{2}U_1.
$$

We then obtain

$$
\nabla^2 Q\vert_W = -J\nabla^W X_Q=\mbox{diag}(1,3/2,\ldots,1).
$$

Note that $V_1$ is normal to $M$ at $p$, as $dQ=\eta \cdot d\eta$ and so $d_pQ(V_1)=1$, so that
$$
\nabla^2 Q\vert_M=\mbox{diag}(1,\ldots,1), 
$$
is the identity, and hence positive definite. This finishes the proof of Proposition \ref{prop:Kepler_problem_convexity}.

\section{Geometric structures on moduli spaces} 

In this section, we construct the foliation of Theorem \ref{thm:IPfoliation}. This is a natural extension of the construction of \cite{M}, and we will refer to that paper for more details, and keep the exposition streamlined for the sake of brevity.

Let $M=\mathbf{OB}(P,\phi)$ be an IP $5$-fold, $\alpha \in \mathbf{Reeb}(P,\phi)$, and $\xi=\ker \alpha$ with Reeb vector field $R_\alpha$. Let $\theta_B: B\backslash L\rightarrow S^1$ be a concrete open book on $B$ adapted to $\alpha_B=\alpha\vert_B$, and $\theta_M: M\backslash B\rightarrow S^1$ a concrete open book adapted to $\alpha$. We denote by $P_\varphi=\overline{\theta_M^{-1}(\varphi)}$, $F_\varphi=\overline{\theta_B^{-1}(\varphi)}$ the $\varphi$-pages. 

\medskip

\textbf{Holomorphic open book over $B$.} As in \cite{M} (see \cite{W10,A11}), the open book at $B$ can be made holomorphic for some $J_B$ compatible with a stable Hamiltonian structure $\mathcal{H}_B=(\lambda_B,d\alpha_B)$ deforming $(\alpha_B,d\alpha_B)$. The result is a finite energy foliation $\mathcal{M}_B$ of $[0,\infty)\times B$, by punctured holomorphic curves which are asymptotically cylindrical Liouville completions $\widehat F_\varphi$ of $F_\varphi$, asymptotic to $L$, so that $\widehat F_\varphi$ projects to $F_\varphi$ under the projection $\theta_B: \mathbb{R}\times B \rightarrow B$. Moreover, we have $\mathcal{M}_B^q:=\mathcal{M}_B/\mathbb R\cong S^1$. All curves are still symplectic with respect to $d\alpha_B$. Similarly as in \cite{M}, using the concrete open book on $M$ and the concrete Lefschetz on $W$, the goal will be to extend this foliation to all of $W$. 

\medskip

\textbf{Holomorphic open book over $M$.} As in \cite{M}, the open book at $M$ can also be made holomorphic for some $J$ compatible with a stable Hamiltonian structure $\mathcal{H}=(\lambda,d\alpha)$ deforming $(\alpha,d\alpha)$. The result is a codimension-2 holomorphic foliation $\mathcal{L}$ of $[0,\infty)\times M$ whose leaves are $\mathcal{L}_B=\mathbb{R}\times B$ and Liouville completions $\widehat P_\varphi$ of the pages $P_\varphi$. We have that $\widehat P_\varphi$ projects to $P_\varphi$ under the projection $p_M: \mathbb R\times M\rightarrow M$, and is asymptotically cylindrical to the $\mathbb{R}$-invariant hypersurface $\mathcal{L}_B$ in the sense of \cite{MS}. Moreover, $J$ can be chosen to coincide with $J_B$ along $\mathcal{L}_B$, and therefore the previously constructed curves are $J$-holomorphic.

\medskip

\textbf{Holomorphic Lefschetz fibration on $W$.} The next step is to extend the foliation $\mathcal{L}$ to the completion $\widehat W$, to a foliation denoted $\widehat{\mathcal L}$. This can be done in such a way that the regular fibers of the concrete Lefschetz fibration $\pi_W: W \rightarrow \mathbb D^2$ (which induces the open book at $M$ by assumption) are all holomorphic and asymptotically cylindrical over the holomorphic hypersurface $\widehat{P}_B=P_B\cup [0,\infty)\times B$, with respect to an extension of $J$ to all of $\widehat W$. This almost complex structure can be further chosen to agree with the standard integrable complex structure near the (non-degenerate) singularities of $\pi_W$. Therefore the foliation $\widehat{\mathcal{L}}$ splits into a union $\widehat{\mathcal{L}}=\widehat{ \mathcal{L}}_{reg}\cup \widehat{\mathcal{L}}_{sing}$ of regular and singular fibers, where the latter correspond to the singular fibers of $\pi_W$.

\medskip

\textbf{Holomorphic foliation of $\widehat W$.} The moduli space $\mathcal{M}_B$ then naturally extends to a moduli space $\mathcal M^W$ on $\widehat W$. We may then consider its SFT-Gromov compactification $\widehat{\mathcal{M}}^W$ obtained by adding strata of nodal curves, as well as its pointed version $\widehat{\mathcal{M}}^W_*$ obtained by adding a single marked point in the domain of each curve, together with the resulting evaluation map $ev: \widehat{\mathcal{M}}^W_* \rightarrow \widehat W$ which evaluates a curve at its marked point, and the forgetful map $\pi_*: \widehat{\mathcal{M}}^W_*\rightarrow \widehat{\mathcal{M}}^W$ which forgets the marked point. As in \cite{M}, Siefring intersection theory (applied iteratively), combined with Wendl's results on fillings of planar open books \cite[Theorem 7, Theorem 8]{W10b}, imply that the image of a curve in $\widehat{\mathcal{M}}^W$ lies completely in a leaf of $\widehat{\mathcal{L}}$, and moreover the curves in $\widehat{\mathcal{M}}^W$ foliate each regular leaf in $\widehat{\mathcal{L}}_{reg}$ as the fibers of a concrete Lefschetz fibration of the given abstract type $\mathbf{LF}(F,\phi_F)$. In particular, there is a subfoliation $\widehat{\mathcal{M}}^B$ foliating $\widehat{W}_B$ and giving such a concrete Lefschetz fibration $\pi^B_*:\widehat{\mathcal{M}}^B_*\rightarrow \widehat{\mathcal{M}}^B\cong \mathbb C$. Exactness of the symplectic form on $W$ implies that theses fibrations are minimal. All curves are immersed, and non-nodal curves are actually embedded.

Then $\widehat{\mathcal{M}}^W_*$ admits a stratification $\widehat{\mathcal{M}}^W_*=\widehat{\mathcal{M}}_*^0\bigsqcup \widehat{\mathcal{M}}_*^1$, where $\widehat{\mathcal{M}}_*^i$ consists of pointed curves with precisely $i$ nodes, and having closure $\bigcup_{j\leq i} \widehat{\mathcal{M}}_*^i$. In particular $\widehat{\mathcal{M}}_*^0$ is the top open stratum consisting of non-nodal curves, and $\widehat{\mathcal{M}}_*^1$ is closed. We similarly have a stratification for the unmarked moduli space $\widehat{\mathcal{M}}^W=\widehat{\mathcal{M}}^0\bigsqcup \widehat{\mathcal{M}}^1$, where $\widehat{\mathcal{M}}^i=\pi_*(\widehat{\mathcal{M}}^i_*)$. The expected dimension of $\widehat{\mathcal{M}}_*^i$ is $6-2i$, and that of $\widehat{\mathcal{M}}^i$ is $4-2i$. An analogous analysis as carried out in \cite[sec.\ 4.7]{Mo} for curves lying in holomorphic hypersurfaces implies that each nodal strata $\widehat{\mathcal{M}}^i$ is Fredholm regular. 

\begin{figure}
    \centering
    \includegraphics[width=0.9\linewidth]{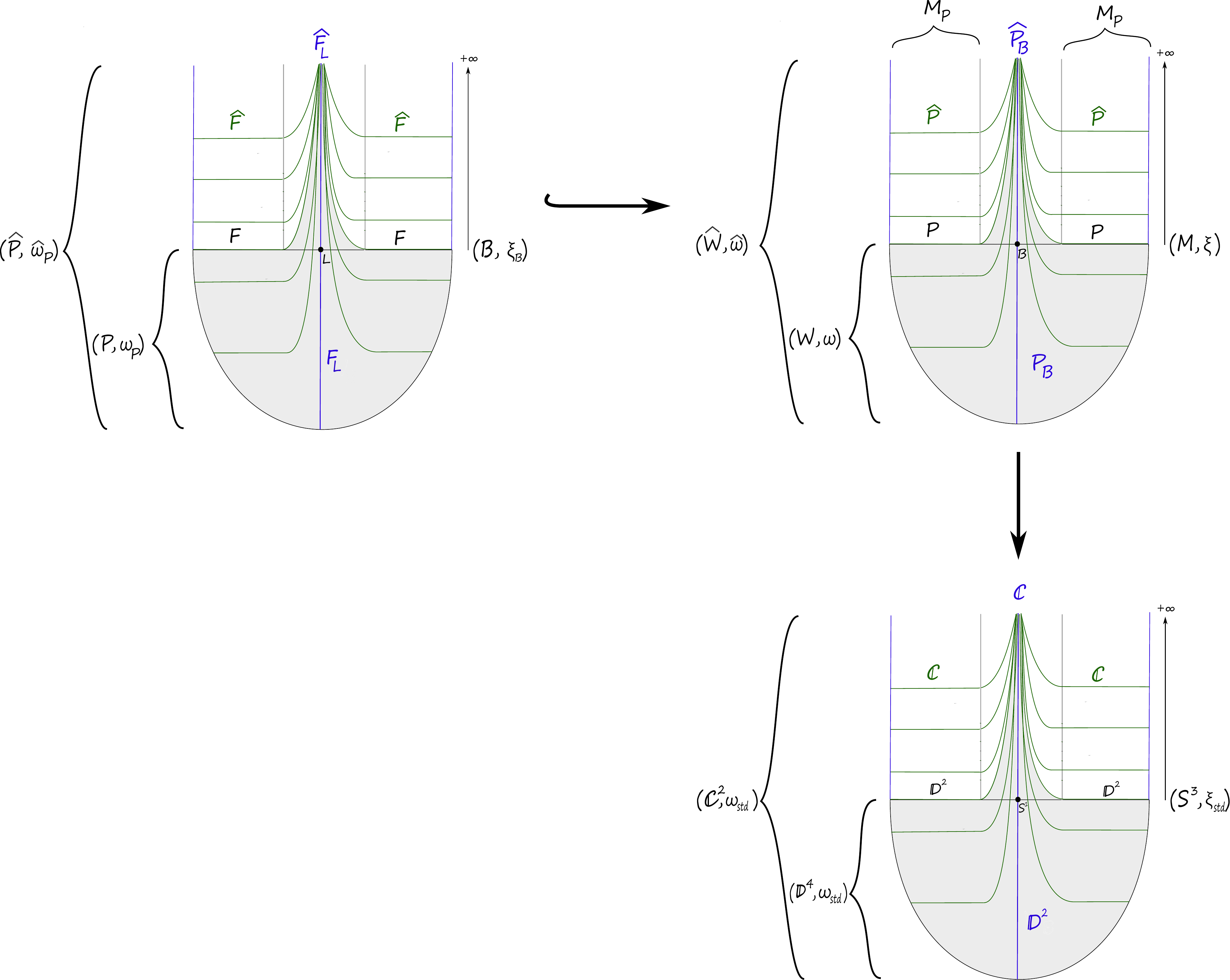}
    \caption{The anatomy of an IP foliation of an IP filling. The non-shaded region is foliated by the curves forming the cylindrical end of $\widehat{\mathcal{M}}$. These curves have constant $\mathbb R$-component over the piece $M_P$ of $M$ lying away from the binding.}
    \label{fig:IP}
\end{figure}

\medskip

\textbf{Structural diagrams.} From the above, we have a diagram of the following form:

\begin{center}
    \begin{tikzcd}
         \widehat P_B \arrow[dr,"\pi^B"]& \arrow[l,"\cong","ev"']\widehat{\mathcal{M}}_*^B \arrow[r,hookrightarrow]\arrow[d,"\pi_*^B"] & \widehat{\mathcal{M}}^W_*\arrow[r, "ev", "\cong"'] \arrow[d, "\pi_*"]&
        \widehat W \arrow[dl, "\pi"]\arrow[r,"\widehat{\pi}_W"]& \mathbb C\\ & \widehat{\mathcal{M}}^B\cong \mathbb C \arrow[r,hookrightarrow]&
        \widehat{\mathcal{M}}^W\cong \mathbb C^2\arrow[r,"\widehat{\pi}_\mathcal{M}"]& \mathbb C &\\
    \end{tikzcd}
\end{center}
The inclusions are inclusions of regular fibers, and all projections are Lefschetz fibrations. See Figure \ref{fig:IP}. There is also a truncation of the above diagrams, correspoding to the filling regions in Figure \ref{fig:IP}, of the form

\begin{center}
    \begin{tikzcd}
         P_B \arrow[dr,"\pi^B"]& \arrow[l,"\cong","ev"']\overline{\mathcal{M}}_*^B \arrow[r,hookrightarrow]\arrow[d,"\pi_*^B"] & \overline{\mathcal{M}}^W_*\arrow[r, "ev", "\cong"'] \arrow[d, "\pi_*"]&
        W \arrow[dl, "\pi"]\arrow[r,"\pi_W"]& \mathbb D^4\\ & \overline{\mathcal{M}}^B\cong \mathbb D^2 \arrow[r,hookrightarrow]&
        \overline{\mathcal{M}}^W\cong \mathbb D^4\arrow[r,"\pi_\mathcal{M}"]& \mathbb D^2 &\\
    \end{tikzcd}
\end{center}

Here, $\overline{\mathcal{M}}_*^B=ev^{-1}(P_B), \overline{\mathcal{M}}_*^W=ev^{-1}(W), \overline{\mathcal{M}}^B=\pi_*(\overline{\mathcal{M}}_*^B), \overline{\mathcal{M}}^W=\pi_*(\overline{\mathcal{M}}_*^W)$, and $\pi_*^B=\pi_*\vert_{\overline{\mathcal{M}}_*^B}$. Note that $\widehat{\mathcal{M}}^W$ has a natural cylindrical end $\widehat{\mathcal{M}}^W\backslash \overline{\mathcal{M}}^W=[0,+\infty)\times \partial \overline{\mathcal{M}}^W$, corresponding those curves completely contained in the cylindrical end $[0,+\infty)\times M$. Here, the boundary $\overline{\mathcal{M}}^q:=\partial \overline{\mathcal{M}}^W\cong S^3$ is identified with the leaf space of the foliation obtained by projecting the latter curves to $M$. And $\overline{\mathcal{M}}^W\cong \mathbb D^4$ is the trivial Lefschetz fibration over $\mathbb D^2$ with fiber $\overline{\mathcal{M}}^B\cong \mathbb D^2$.

\medskip

This finishes the proof of Theorem \ref{thm:IPfoliation}.

\medskip

\textbf{Symplectic structure.} We now construct a symplectic form on $\overline{\mathcal{M}}^W$, and prove Theorem \ref{thm:symplectic_structure}, i.e.\ that it gives a symplectic filling of its boundary. We follow the same approach as in \cite{M}. Let $$\mathcal{P}=\{\varphi \in C^\infty((0,+\infty),(0,1)): \varphi^\prime>0\}$$ be the space of orientation-preserving diffeomorphisms between $(0,+\infty)$ and $(0,1)$. For $\varphi \in \mathcal{P},$ we let
$$
\omega^\varphi:=\left\{\begin{array}{cc}
   \omega  & \mbox{ on } W \\
   d(e^{\varphi(t)}\alpha)  & \mbox{ on } (0,+\infty)\times M.
\end{array}\right.
$$
This is simply a truncation of the symplectic form $\widehat{\omega}$ on $\widehat W$. Define $\omega_*^\varphi=ev^*\omega^\varphi \in \Omega^2(\widehat{\mathcal{M}}_*)$. If now $u\in \widehat{\mathcal{M}}$ and $v,w\in T_u\widehat{\mathcal{M}}$, we define
$$
(\omega^\varphi_{\mathcal{M}})_{u}(v,w)=\int_{z\in F_u} (\omega_*^\varphi)_{u(z)}(v(z),w(z))dz,
$$
where we denote $F_u=\pi_*^{-1}(u)$ the domain of $u$, and $dz=\omega_*^\varphi\vert_{F_u}$. This is also well-defined over the nodal strata. Then  $\omega^\varphi_{\mathcal{M}}$ is symplectic, as follows from the geometry of the foliation (i.e.\ the fact that the symplectic normal bundles to the curves are positive), and of the form
$$
\omega^\varphi_{\mathcal{M}}=\left\{\begin{array}{cc}
   \omega_{\mathcal{M}}  & \mbox{ on } \overline{\mathcal{M}}^W\\
   d(e^{\varphi(t)}\alpha_\mathcal{M})  & \mbox{ on } [0,+\infty)\times \partial \overline{\mathcal{M}}^W,
\end{array}\right.
$$
where $(\alpha_\mathcal{M})_u(v)=\int_{z \in F_u}\alpha_{u(z)}(v(z))dz$ is a contact form on $\partial \overline{\mathcal{M}}^W\cong S^3$, the \emph{holomorphic shadow} of $\alpha$. We refer to \cite{M} for more properties of this contact form.

\medskip

This finishes the proof of Theorem \ref{thm:symplectic_structure}.

\medskip

\textbf{Convexity.} We now show that the holomorphic shadow $(\partial \overline{\mathcal{M}}^W,\alpha_\mathcal{M})$ is strictly convex in $(\overline{\mathcal{M}}^W,\omega_\mathcal{M})$, if $(M,\alpha)$ is strictly convex in $(W,\omega)$. Let $H: N \rightarrow \mathbb R$ be a Hamiltonian as in Definition \ref{def:convexity}, defined on a collar neighbourhood $N=(-\epsilon,\epsilon)\times M$ of $M$ in $W$. Let $\mathcal{N}_*=ev^{-1}((-\epsilon,\epsilon)\times M_P)\subset \widehat{\mathcal{M}}_*^W$, where $M_P\subset M$ corresponds to the mapping torus part of the open book, and where the holomorphic curves can be taken by construction to lie in the slices $\{a\}\times M_P$; see Figure \ref{fig:IP}.

and $\mathcal{N}=\pi_*(\mathcal{N}_*)\subset \widehat{\mathcal{M}}^W$. Define
$$
H_\mathcal{M}: \mathcal{N}\rightarrow \mathbb R,
$$
$$ H_\mathcal{M}(u)=\int_{z\in F_u \cap \mathcal{N}_*} H(u(z))dz.
$$
Note that, by choice of $\mathcal{N}_*$, we have $H_\mathcal{M}^{-1}(0)=\partial \overline{\mathcal{M}}^W$. We take the metric $g=\omega(\cdot, J^\xi\cdot)$ on $\widehat W$, where $J^\xi$ is an $\alpha$-compatible almost complex structure canonically associated with the SHS $J$ used for the foliation (see \cite{M}). We let $\nabla$ be the associated Levi-Civita connection. This induces a metric $g_\mathcal{M}$ on $\widehat{\mathcal{M}}$, given by
$$
(g_\mathcal{M})_u(v,w)=\int_{z\in F_u} g_{u(z)}(v(z),w(z))dz,
$$
for $v,w \in T_u\widehat{\mathcal{M}}$. We let $\nabla_\mathcal{M}$ be the associated Levi-Civita connection. Then the Hessian of $H_\mathcal{M}$ is given by
$$
\nabla_\mathcal{M}^2H_\mathcal{M}(v,w)=\int_{z\in F_u\cap \mathcal{N}_*} (\nabla^2H)_{u(z)}(v(z),w(z))dz.
$$
The portion of the curves in $H_\mathcal{M}^{-1}(0)$ that are mapped to $(-\epsilon,\epsilon)\times M_P$ actually lie in $\{0\}\times M_P\subset M$, and so along $H_\mathcal{M}^{-1}(0)$ we have $\nabla_\mathcal{M}^2H_\mathcal{M}>0$, since $\nabla^2H>0$ pointwise on $M$.

\medskip

This finishes the proof of Theorem \ref{thm:convexity}.

\end{document}